\documentclass[a4paper,english,11pt]{amsart}
\usepackage{babel}
\usepackage{microtype,txfonts}
\usepackage{hyperref}
\newtheorem{theorem}{Theorem}[section]

\newtheorem{example}{\it Example\/}

\def\F{\mathbf{F}}\def\P{\mathbf{P}}\def\G{\mathbf{G}}
\def\rk{\mathrm{rk}}
\let\og\guillemotleft\let\fg\guillemotright
\title{Flag bundles, Segre polynomials and push-forwards.}
\author{Lionel Darondeau}
\email{lionel.darondeau@normalesup.org}
\author[Piotr Pragacz]{Piotr Pragacz$^\dagger$}
\email{P.Pragacz@impan.pl}
\address{Institut de Math\'ematiques de l'Acad\'emie des Sciences de Pologne\\
ul. \'Sniadeckich 8, 00-656 Warszawa --- Polska}
\thanks{$^\dagger$Supported by National Science Center (NCN) grant no. 2014/13/B/ST1/00133}
\begin{document}
\begin{abstract}
  In this note, we give Gysin formulas for partial flag bundles for the classical groups.
  We then give Gysin formulas for Schubert varieties in Grassmann bundles, including isotropic ones.
  All these formulas are proved in a rather uniform way by using the step-by-step construction of flag bundles and the Gysin formula for a projective bundle.
  In this way we obtain a comprehensive list of new universal formulas.
  The content of this paper was presented by Piotr Pragacz at the International Festival in Schubert Calculus in Guangzhou, November 6-10, 2017.
\end{abstract}
\maketitle
\section{Introduction}
\label{se:intro}

Let \(E\to X\) be a vector bundle of rank \(n\) on a variety \(X\) over an algebraically closed field.
Let \(\pi\colon\F(E)\to X\) be the bundle of flags of subspaces of dimensions \(1,2,\dots, n-1\) in the fibers of \(E\to X\).
The flag bundle \(\F(E)\) is used, \textit{e.g.}, in splitting principle, a standard technique which reduces questions about vector bundles to the case of line bundles; namely the pullback bundle \(\pi^{\ast}E\) decomposes as a direct sum of line bundles.
One can construct \(\F(E)\) inductively as a sequence of projective bundles, using the following iterative step, that decreases the rank by \(1\).
Let \(p_{1}\colon\P(E)\to X\) denote the projective bundle of lines in the fibers of \(E\), and let \(U_{1}:=\mathcal{O}_{\P(E)}(-1)\) denote the universal subbundle on \(\P(E)\), then one has a universal exact sequence
of vector bundles on \(\P(E)\)
\[
0\to U_{1}\longrightarrow p_{1}^{\ast}E \longrightarrow Q_{n-1}\to 0,
\]
where \(Q_{n-1}\) (the universal quotient bundle on \(\P(E)\)) is a rank \(n-1\) vector bundle.
Replacing \(E\) by \(Q_{n-1}\), one obtains a universal subbundle on \(\P(Q_{n-1})\), that it is convenient to denote \(U_{2}/U_{1}\), together with a universal quotient bundle \(Q_{n-2}\). Iterating this process until obtaining a quotient bundle \(Q_{1}\) of rank one, one gets a sequence of projective bundles
\begin{equation}
  \label{eq:full}
  \F(E)
  :=
  \P(Q_{2})
  \stackrel{p_{n-1}}\longrightarrow
  \cdots
  \to
  \P(Q_{n-1})
  \stackrel{p_{2}}\longrightarrow
  \P(E)
  \stackrel{p_{1}}\longrightarrow
  X,
\end{equation}
and universal exact sequences of vector bundles on \(\P(Q_{n-i+1})\):
\begin{equation}
  \label{eq:taut}
  0
  \to
  U_{i}/U_{i-1}
  \to
  p_{i}^{\ast}Q_{n-i+1}
  \to
  Q_{n-i}
  \to
  0
\end{equation}
such that
\(
\pi^{\ast}E
=
U_{1}\oplus U_{2}/U_{1}\oplus\cdots\oplus U_{n-1}/U_{n-2}\oplus Q_{1}
\).
Here \(U_{i+1}/U_{i}\) stands for the universal subbundle on \(\P(Q_{n-i+1})\) and also for the pullbacks of these to \(\F(E)\).

Now we would like to outline how to obtain a Gysin formula for the flag bundle \(\pi\colon\F(E)\to X\) (\textit{cf.} Example \ref{ex1}),
and introduce some notation.

We shall work in the framework of intersection theory of \cite{Fulton}.
Recall that a proper morphism \(g\colon Y\to X\) of nonsingular algebraic varieties over an algebraically closed field yields an additive map \(g_{\ast}\colon A^{\bullet}Y\to A^{\bullet}X\) of Chow groups induced by push-forward cycles, called the Gysin map.
The theory developed in \cite{Fulton} allows also one to work with singular varieties, or with cohomology.
In this note, \(X\) will always be nonsingular.

For \(E\to X\) a vector bundle, let \(s(E)\) be the Segre class of \(E\), that is the formal inverse of the Chern class \(c(E)\) in the Chow ring of \(X\).
Let \(\xi=c_{1}(\mathcal{O}_{\P(E)}(1))\); then \(A^{\bullet}(\P(E))\) is generated algebraically by \(\xi\) over \(A^{\bullet}X\)---here we identify \(A^{\bullet}X\) with a subring of \(A^{\bullet}(\P(E))\)---and
\begin{equation}
  \label{eq:segre}
  (p_{1})_{\ast} \xi^{i}
  =
  s_{i-(n-1)}(E),
\end{equation}
\textit{cf.} \cite{Fulton}.
To obtain a Gysin formula for the sequence of projective bundles (\ref{eq:full}), it suffices to appropriately iterate formula (\ref{eq:segre}).
The intermediate formulas involve the individual Segre classes of the universal quotient bundles, that can be eliminated using (\ref{eq:taut}) and the Whitney sum formula.
However, it seems rather difficult to obtain a universal formula in this way.
A universal formula should hold for \emph{any} polynomial in characteristic classes of universal vector bundles and depend explicitly on the Segre classes of the original bundle \(E\).
To obtain such a formula, we use the generating series of the Segre classes of the universal quotient bundles.
A prototype is the reformulation of (\ref{eq:segre}) in
\begin{equation}
  \label{eq:fund}
  (p_{1})_{\ast} \xi^{i}
  =
  [t^{n-1}]
  \big(t^{i}s_{1/t}(E)\big),
\end{equation}
where we consider the specialization in \(x=1/t\) of the Segre polynomial \(s_{x}(E)=\sum_{i}s_{i}(E)x^{i}\) and where for a monomial \(m\) and a Laurent polynomial \(P\), \([m](P)\) denotes the coefficient of \(m\) in \(P\).
Formula (\ref{eq:fund}) and the projection formula imply that for any polynomial \(f\) in one variable with coefficients in \(A^{\bullet}X\)
\begin{equation}
  \label{eq:fund_f}
  (p_{1})_{\ast} f(\xi)
  =
  [t^{n-1}]
  \big(f(t)s_{1/t}(E)\big).
\end{equation}
In this formula, (i) one does not need to expand \(f\) into a combination of monomials; (ii) one uses the Segre polynomial that, like the total Segre class, is a group homomorphism from the Grothendieck group of \(X\) to the multiplicative group of units with degree zero term \(=1\) in \(A^{\bullet}X\).

Iterating the Gysin formula (\ref{eq:fund_f}) yields a closed universal Gysin formula for the flag bundle \(\F(E)\to X\), as announced in Example \ref{ex1}.

It is clear that the outlined strategy of proof applies to more general step-by-step constructions than the construction (\ref{eq:full}) of the flag bundle \(\F(E)\to X\).
Considering the truncated composition \(p_{k}\circ\cdots \circ p_{1}\) in (\ref{eq:full}) yields formulas for
full flag bundles, \textit{i.e.} bundles of flags of subspaces of dimensions \(1,2,\dots,k\) in the fibers, for \(k=1,\dots,n-1\).
Then, using certain commutative diagrams (see ~\cite[(5)]{DP1}), one extends these formulas to arbitrary partial flag bundles.

One other interesting generalization is to restrict to the zero locus of a section of some vector bundle
at each step of the sequence of projective bundles.
In other words, one can impose some geometric conditions that the subspaces of the flag have to satisfy.
An illustrative example is Theorem~\ref{thm:gysin-BD}, in the orthogonal setting, obtained by considering at each step quadric bundles of isotropic lines in projective bundles of lines.

This method of step-by-step construction of generalized flag bundles leads to uniform short proofs of the different results announced in this note.

This note is organized as follows.
In Section \ref{se:flags}, we shall announce universal Gysin formulas for partial flag bundles for general linear groups, symplectic groups and orthogonal groups.
The proofs of the results announced there can be found in \cite{DP1}.

In Section~\ref{se:schubert} we give Gysin formulas for Kempf-Laksov flag bundles.
These generalized flag bundles are used to desingularize Schubert varieties in Grassmann bundles.
Theorem~\ref{thm:gysin-KL-A} is established in~\cite{DP2}.
Theorem~\ref{thm:gysin-KL-C} is announced for the first time in the present note.

\section{Universal Gysin formulas for flag bundles}
\label{se:flags}

In this section, the letter \(f\) denotes a polynomial in the indicated number of variables with coefficients in \(A^{\bullet}X\).
The appropriate symmetries that \(f\) has to satisfy to be in the Chow ring of the flag bundle under consideration are always implied.

We shall discuss separately the cases of general linear groups, symplectic groups and orthogonal groups.

\subsection{General linear groups}

Let \(E\to X\) be a rank \(n\) vector bundle.
Let \(1\leq d_{1}<\cdots< d_{m}=d\leq n-1\) be a sequence of integers.
We denote by \(\pi\colon\F(d_{1},\dots,d_{m})(E)\to X\)
the bundle of flags of subspaces of dimensions \(d_{1},\dots,d_{m}\) in the fibres of \(E\).
On \(\F(d_{1},\dots,d_{m})(E)\), there is a universal flag
\(U_{d_{1}}\subsetneq\cdots\subsetneq U_{d_{m}}\)
of subbundles of \(\pi^{\ast}E\), where \(\rk(U_{d_{k}})=d_{k}\)
(the fiber of \(U_{d_{k}}\) over the point
\((V_{d_{1}}\subsetneq\cdots\subsetneq V_{d_{m}}\subset E_{x})\),
where \(x\in X\), is equal to \(V_{d_{k}}\)).
For a foundational account on flag bundles, see~\cite{Groth2}.

For \(i=1,\dots,d\), set \(\xi_{i}=-c_{1}(U_{d+1-i}/U_{d-i})\).
\medskip
\begin{theorem}
  \label{thm:gysin-A}
  With the above notation,
  for
  \(
  f(\xi_{1},\dots,\xi_{d})
  \in
  A^{\bullet}(\F(d_{1},\dots,d_{m})(E)),
  \)
  one has
  \[
  \pi_{\ast}f(\xi_{1},\dots,\xi_{d})
  =
  \Big[
  {t_{1}}^{e_{1}}\dots{t_{d}}^{e_{d}}
  \Big]
  \bigg(
  {
  \textstyle
  f(t_{1},\dots,t_{d})\,
  \prod\limits_{1\leq i<j\leq d}
  (t_{i}-t_{j})
  \prod\limits_{1\leq i \leq d}
  s_{1/t_{i}}(E)
  }
  \bigg),
  \]
  where for \(j=d-d_{k}+i\) with \(i=1,\dots,d_{k}-d_{k-1}\),
  we denote
  \(
  e_{j}
  =
  n-i
  \).
\end{theorem}
\medskip
\begin{example}
  \label{ex1}
  For the complete flag bundle \(\pi\colon\F(E)\to X\), one has
  \[
  \pi_{\ast}f(\xi_{1},\dots,\xi_{n-1})
  =
  \Big[{\textstyle\prod\limits_{i=1}^{n-1} t_{i}^{(n-1)}}\Big]
  \bigg(
  {\textstyle
  f(t_{1},\dots,t_{n-1})
  \prod\limits_{1\leq i< j\leq n-1}
  (t_{i}-t_{j})
  \prod\limits_{i=1}^{n-1}
  s_{1/t_{i}}(E)
  }
  \bigg);
  \]
  and for the Grassmann bundle \(\pi\colon\F(d)(E)\to X\), one has
  \[
  \pi_{\ast}f(\xi_{1},\dots,\xi_{d})
  =
  \Big[{\textstyle\prod\limits_{i=1}^{d}t_{i}^{n-i}}\Big]
  \bigg(
  {\textstyle
  f(t_{1},\dots,t_{d})
  \prod\limits_{1\leq i< j\leq d}
  (t_{i}-t_{j})
  \prod\limits_{i=1}^{d}
  s_{1/t_{i}}(E)
  }
  \bigg).
  \]
\end{example}

\subsection{Symplectic groups}
Let \(E\to X\) be a rank \(2n\) vector bundle equipped with a non-degenerate symplectic form \(\omega\colon E\otimes E\to L\) (with values in a certain line bundle \(L\to X\)). We say that a subbundle \(S\) of \(E\) is isotropic if \(S\) is a subbundle of its symplectic complement \(S^{\omega}\), where
\[
S^{\omega}
:=
\{w\in E\mid \forall v\in S\colon \omega(w,v)=0\}.
\]

Let \(1\leq d_{1}<\cdots<d_{m}\leq n\) be a sequence of integers.
We denote by \(\pi\colon\F^{\omega}(d_{1},\dots,d_{m})(E)\to X\) the bundle of flags of isotropic subspaces of dimensions \(d_{1},\dots,d_{m}\) in the fibers of \(E\).
On \(\F^{\omega}(d_{1},\dots,d_{m})(E)\), there is a universal flag \(U_{d_{1}}\subsetneq\cdots\subsetneq U_{d_{m}}\) of subbundles of \(\pi^{\ast}E\), where \(\rk(U_{d_{k}})=d_{k}\).

For \(i=1,\dots,d\), set \(\xi_{i}=-c_{1}(U_{d+1-i}/U_{d-i})\).
\medskip
\begin{theorem}
  \label{thm:gysin-C}
  With the above notation,
  for
  \(
  f(\xi_{1},\dots,\xi_{d})
  \in
  A^{\bullet}(\F^{\omega}(d_{1},\dots,d_{m})(E)),
  \)
  one has
  \[
  \pi_{\ast}f(\xi_{1},\dots,\xi_{d})
  =
  \Big[
  {t_{1}}^{e_{1}}\cdots{t_{d}}^{e_{d}}
  \Big]
  \bigg(
  {
  \textstyle
  f(t_{1},\dots,t_{d})\,
  \prod\limits_{1\leq i<j\leq d}
  (c_{1}(L)+t_{i}+t_{j})
  (t_{i}-t_{j})
  \prod\limits_{1\leq i \leq d}
  s_{1/t_{i}}(E)
  }
  \bigg),
  \]
  where for \(j=d-d_{k}+i\) with \(i=1,\dots,d_{k}-d_{k-1}\),
  we denote
  \(
  e_{j}
  =
  2n-i
  \).
\end{theorem}
\medskip

\begin{example}
  For the symplectic Grassmann bundle \(\pi\colon\F^{\omega}(d)(E)\to X\), where \(\omega\) has values in a trivial line bundle, one has
  \[
  \pi_{\ast}f(\xi_{1},\dots,\xi_{d})
  =
  \Big[{\textstyle\prod\limits_{i=1}^{d}t_{i}^{2n-i}}\Big]
  \bigg(
  {\textstyle
  f(t_{1},\dots,t_{d})
  \prod\limits_{1\leq i< j\leq d}
  (t_{i}^{2}-t_{j}^{2})
  \prod\limits_{i=1}^{d}
  s_{1/t_{i}}(E)
  }
  \bigg).
  \]
\end{example}

\subsection{Orthogonal groups}
Let \(E\to X\) be a vector bundle of rank \(2n\) or \(2n+1\) equipped with a non-degenerate orthogonal form \(Q\colon E\otimes E\to L\) (with values in a certain line bundle \(L\to X\)).  We say that a subbundle \(S\) of \(E\) is isotropic if \(S\) is a subbundle of its orthogonal complement \(S^{\perp}\), where
\[
S^{\perp}
:=
\{w\in E\mid \forall v\in S\colon Q(w,v)=0\}.
\]

Let \(1\leq d_{1}<\cdots<d_{m}\leq n\) be a sequence of integers. We denote by \(\pi\colon \F^{Q}(d_{1},\dots,d_{m})(E)\to X\) the bundle of flags of isotropic subspaces of dimensions \(d_{1},\dots,d_{m}\) in the fibers of \(E\). On \(\F^{Q}(d_{1},\dots,d_{m})(E)\), there is a universal flag \(U_{d_{1}}\subsetneq\dots\subsetneq U_{d_{m}}\) of subbundles of \(\pi^{\ast}E\), where \(\rk(U_{d_{k}})=d_{k}\).

For \(i=1,\dots,d\), set \(\xi_{i}=-c_{1}(U_{d+1-i}/U_{d-i})\).
\medskip
\begin{theorem}
  \label{thm:gysin-BD}
  With the above notation,
  for
  \(
  f(\xi_{1},\dots,\xi_{d})
  \in
  A^{\bullet}(\F^{Q}(d_{1},\dots,d_{m})(E)),
  \)
  one has
  \begin{eqnarray*}
    &&\pi_{\ast}f(\xi_{1},\dots,\xi_{d})
    = \\ &&\qquad
    \Big[
    {t_{1}}^{e_{1}}\cdots{t_{d}}^{e_{d}}
    \Big]
    \bigg(
    {
    \textstyle
    f(t_{1},\dots,t_{d})\,
    \prod\limits_{1\leq i \leq d}
    (2t_{i}+c_{1}(L))
    \prod\limits_{1\leq i<j\leq d}
    (c_{1}(L)+t_{i}+t_{j})
    (t_{i}-t_{j})
    \prod\limits_{1\leq i \leq d}
    s_{1/t_{i}}(E)
    }
    \bigg),
  \end{eqnarray*}
  where for \(j=d-d_{k}+i\) with \(i=1,\dots,d_{k}-d_{k-1}\),
  we denote
  \(
  e_{j}
  =
  \rk(E)-i
  \).
\end{theorem}
\medskip

Note that, if the rank is \(2n\) and \(d=n\), we consider \emph{both} of the two isomorphic connected components of the flag bundle. Thus, if one is interested in only one of the two components, the result should be divided by \(2\). When \(c_{1}(L)=0\), this makes appear the usual coefficient \(2^{n-1}\).

\section{Universal Gysin formulas for Kempf--Laksov flag bundles}
\label{se:schubert}
In this section, we give Gysin formulas for Kempf--Laksov flag bundles, that are desingularizations of Schubert bundles in Grassmann bundles. We also extend the results to the symplectic setting. The orthogonal cases will be treated elsewhere.

\subsection{General linear groups}

Let \(E\to X\) be a rank \(n\) vector bundle on a variety \(X\) with a reference flag of bundles
\(E_{1}\subsetneq\cdots\subsetneq E_{n}=E\) on it, where \(\rk(E_{i})=i\).
Let \(\pi\colon\G_{d}(E)=\F(d)(E)\to X\) be the Grassmann bundle of subspaces of dimension \(d\) in the fibers of \(E\).
For any partition \(\lambda\subseteq(n-d)^{d}\), one defines the \textsl{Schubert bundle} \(\varpi_{\lambda}\colon\Omega_{\lambda}(E_{\bullet})\to X\) in \(\G_{d}(E)\) over the point \(x\in X\)
by
\begin{equation}
  \label{eq:def_omega}
  \Omega_{\lambda}(E_{\bullet})(x)
  :=
  \{
  V\in \G_{d}(E)(x)
  \colon
  \dim(V\cap E_{n-d-\lambda_{i}+i}(x))\geq i,
  \mbox{ for }i=1,\dots,d
  \}.
\end{equation}

We denote by
\[
(\nu_{1},\dots,\nu_{d})
:=
(n-d-\lambda_{d}+d,\dots,n-d-\lambda_{1}+1)
\]
the dimensions of the spaces of the reference flag involved in the definition of \(\Omega_{\lambda}(E_{\bullet})\)---in reverse order---.
The partition \(\nu\) is a strict partition, and furthermore, \(n-i\leq\nu_{i}\leq\nu_{1}=n-\lambda_{d}\leq n\) for any \(i\).
Note that the above definition of \(\Omega_{\lambda}(E_{\bullet})\) can be restated using \(\nu\) with the conditions
\begin{equation}
  \label{eq:condition_nu}
  \dim(V\cap E_{\nu_{i}}(x))\geq d+1-i,\mbox{ for }i=1,\dots,d.
\end{equation}

For a strict partition \(\mu\subseteq(n)^d\) with \(d\) parts,
consider the flag bundle \(\vartheta_{\mu}\colon F_{\mu}(E_{\bullet})\to X\) defined over the point \(x\in X\) by
\begin{equation}
  \label{eq:def_f}
  F_{\mu}(E_{\bullet})(x)
  :=
  \Big\{
  0\subsetneq V_{1}\subsetneq\cdots\subsetneq V_{d} \in \F(1,\dots,d)(E)(x)
  \colon
  V_{d+1-i}\subseteq E_{\mu_{i}}(x),
  \mbox{ for }i=1,\dots,d
  \Big\}.
\end{equation}
We will call \textsl{Kempf--Laksov flag bundles} such bundles \(\vartheta_{\mu}\) introduced in \cite{KL}.

These appear naturally as desingularizations of Schubert bundles.
For a partition \(\lambda\subseteq(n-d)^{d}\), defining \(\nu\) as above, by (\ref{eq:condition_nu}) the forgetful map \(\F(1,\dots,d)(E)\to\G_{d}(E)\) induces a birational morphism \(F_{\nu}(E_{\bullet})\to\Omega_{\lambda}(E_{\bullet})\). On the \textsl{Schubert cell} defined over the point \(x\in X\) by
\[
\mathring\Omega_{\lambda}(E_{\bullet})(x)
:=
\Big\{
V\in \G_{d}(E)(x)
\colon
\dim(V\cap E_{\nu_{i}}(x))= d+1-i,\mbox{ for } i=1,\dots,d
\Big\},
\]
which is open dense in \(\Omega_{\lambda}(E_{\bullet})\), the inverse map is
\(
V\mapsto (V\cap E_{\nu_{d}}(x),\dots, V\cap E_{\nu_{1}}(x))
\).
It establishes a desingularization of \(\Omega_{\lambda}(E_{\bullet})\) (see~\cite{KL}).

Consider the sequence of projective bundles
\begin{equation}\label{sequence}
  F_\mu(E_\bullet)=\P(E_{\mu_1}/U_{d-1}) \to \P(E_{\mu_2}/U_{d-2}) \to \cdots \to \P(E_{\mu_{d-1}}/U_{1}) \to \P(E_{\mu_d})\,,
\end{equation}
where \(U_{d-i+1}/U_{d-i}\) is the universal line bundle on \(\P(E_{\mu_i}/U_{d-i})\). Set \(\xi_i=-c_1(U_{d-i+1}/U_{d-i})\), \(i=1,\dots,d\).

\smallskip
Let \(f\) be a polynomial in \(d\) variables with coefficients in \(A^{\bullet}(X)\).

\medskip
\begin{theorem}
  \label{thm:gysin-KL-A}
  With the above notation, one has
  \[
  (\vartheta_{\mu})_{\ast} f(\xi_1,\ldots,\xi_d)
  =
  \Big[t_{1}^{\mu_{1}-1}\cdots t_{d}^{\mu_{d}-1}\Big]
  \left(
  {\textstyle
  f(t_{1},\dots,t_{d})
  \prod\limits_{1\leq i<j\leq d}(t_{i}-t_{j})
  \prod\limits_{1\leq i\leq d}s_{1/t_{i}}(E_{\mu_{i}})
  }
  \right).
  \]
\end{theorem}
\medskip
A proof of this theorem is based on (\ref{sequence}) and (\ref{eq:fund_f}).

\subsection{Symplectic groups}

Let \(E\to X\) be a rank \(2n\) symplectic vector bundle endowed with the symplectic form \(\omega\colon E\otimes E\to L\) with value in a line bundle \(L\to X\), over a variety \(X\).
For \(d\in\{1,\dots,n\}\), let \(\G_{d}^{\omega}(E)=\F^{\omega}(d)(E)\) be the Grassmannian bundle of isotropic \(d\)-planes in the fibers of \(E\).
Let
\[
0=E_{0}\subsetneq E_{1}\subsetneq\cdots\subsetneq E_{n}=E_{n}^{\omega}\subsetneq\cdots\subsetneq E_{0}^{\omega}=E
\]
be a reference flag of isotropic subbundles of \(E\) and their symplectic complements, where \(\rk(E_{i})=i\). For \(i=1,\dots,n\), we set \(E_{n+i}:= E_{n-i}^{\omega}\).
For a partition \(\lambda\subseteq(2n-d)^{d}\), one defines the Schubert cell \(\mathring{\Omega}_{\lambda}(E_{\bullet})\) in \(\G_{d}^{\omega}(E)\) over the point \(x\in X\) by the conditions
\[
\mathring{\Omega}_{\lambda}(E_{\bullet})(x)
:=
\big\{
V\in\G_{d}^{\omega}(E)(x)
\colon
\dim\big(V\cap E_{2n-d+i-\lambda_{i}}(x)\big)=i,
\mbox{ for }i=1,\dots,d
\big\}.
\]

Denote \(\nu_{d+1-i}:=2n-d+i-\lambda_{i}\) the dimension of the reference space appearing in the \(i\)th condition.
A partition indexing the Schubert cell \(\mathring\Omega_{\lambda}\) must satisfy the conditions \(\nu_{i}+\nu_{j}\neq 2n+1\)
(see \cite[p. 174]{P}, where this is shown for \(d=n\), and for arbitrary \(d\) the argument is the same).
For such partitions one defines the Schubert bundle \(\varpi_{\lambda}\colon \Omega_{\lambda}\to X\) as the Zariski-closure of \(\mathring{\Omega}_{\lambda}\), given over a point \(x\in X\) by the conditions
\[
\Omega_{\lambda}(E_{\bullet})(x)
:=
\big\{
V\in\G_{d}^{\omega}(E)(x)
\colon
\dim\big(V\cap E_{2n-d+i-\lambda_{i}}(x)\big)\geq i,
\mbox{ for }i=1,\dots,d
\big\}.
\]

For a strict partition \(\mu\subseteq(2n)^{d}\) with \(d\) parts, such that \(\mu_{i}+\mu_{j}\neq 2n+1\) for all \(i,j\), we define the \textsl{isotropic Kempf--Laksov bundle} \(\vartheta_{\mu}\colon F_{\mu}(E_{\bullet})\to X\) over the point \(x\in X\) by
\[
F_{\mu}(E_{\bullet})(x)
:=
\big\{
0\subsetneq V_{1}\subsetneq\dots\subsetneq V_{d}\in \F^{\omega}(1,\dots,d)(E)(x)
\colon
V_{d+1-i}\subseteq E_{\mu_{i}}(x)
\big\}.
\]
Note that as in the previous section, \(F_{\nu}(E_{\bullet})\) is birational to \(\Omega_{\lambda}(E_{\bullet})\), but here it is not smooth in general.

Let \(U_{i}\) stands for the restriction to \(F_{\mu}(E_{\bullet})\) of the rank \(i\) universal bundle on \(\F(1,\ldots,d)(E)\).
Set \(\xi_{i}=-c_{1}(U_{d-i+1}/U_{d-i})\), for \(i=1,\ldots,d\).

\smallskip

Let \(f\) be a polynomial in \(d\) variables with coefficients in \(A^{\bullet}(X)\).

\medskip
\begin{theorem}
  \label{thm:gysin-KL-C}
  With the above notation, one has
  \begin{eqnarray*}
    &&
    (\vartheta_{\mu})_{\ast} f(\xi_1,\dots,\xi_d)
    =
    \\
    &&
    \hskip15mm
    \big[t_{1}^{\mu_{1}-1}\cdots t_{d}^{\mu_{d}-1}\big]
    \Big(
    {\textstyle
    f(t_{1},\dots,t_{d})
    \prod\limits_{1\leq i<j\leq d}\!
    (t_{i}-t_{j})
    \prod\limits_{{\scriptstyle 1\leq i<j\leq d\atop\scriptstyle\mu_{i}+\mu_{j}>2n+1}}\hskip-3.5mm
    (c_{1}(L)+t_{i}+t_{j})
    \prod\limits_{1\leq j\leq d}\!
    s_{1/t_{j}}(E_{\mu_{j}})
    }
    \Big).
  \end{eqnarray*}
\end{theorem}
\medskip
A proof of this theorem will appear in a separate publication.

\vfill
\end{document}